\documentclass[12pt,twoside]{amsart}
\usepackage{amsmath}
\usepackage{amssymb}
\usepackage{amscd}
\usepackage{xypic}
\xyoption{all}
\title[Vertex algebras and Hodge structures, (Beilinson-Bernstein correspondence)]{Vertex algebras and Hodge structures, (Beilinson-Bernstein correspondence)}
\author{Mohammad Reza Rahmati \\ \\ \\ \\ \\ }
\thanks{}
\address{Universudad deLaSalle Bajio, Leon, Mexico
\hfill\break 
\hfill\break \\
\hfill\break }
\email{mrahmati@cimat.mx}

\newcommand{\comments}[1]{}

\newtheorem{theorem}{Theorem}[section]

\newtheorem{definition}[theorem]{Definition}
\newtheorem{remark}[theorem]{Remark}

\keywords{Vertex operator algebras, Hodge structure, Highest weight modules, $D$-modules, KZ-equation, Affine Lie algebras, Root systems, Beilinson-Bernstein correspondence, Geometric Langlands correspondence, Wakimoto modules, $\mathfrak{g}$-oper, Virasoro algebra, Fock modules}

\subjclass{}

\begin{document}

\begin{abstract}
We compare the context of Hodge structures with that of vertex algebras of conformal field theory. Vertex algebras appear as the highest weight representations of infinite dimensional Lie algebras. A correspondence between Higgs bundles and opers already is known as non-abelian Hodge theorem due to C. Simpson. The Beilinson-Bernstein localization (correspondence) also compares the context of variation of Hodge structure with that of highest weight modules over flag manifolds of semisimple Lie groups. A more general analogue of the Bernstein correspondence over a local manifold can also be formulted in the context of geometric Langlands correspondence. We discuss a generalized version of Harish-Chandra modules called Wakimoto modules and a generalized Harish-Chandra homomorphism. This text is mainly an expository discussion with a new insight toward the two concepts. We conclude with an explanation of geometric Langlands correspondence. 
\end{abstract}

\maketitle

%\noindent Version: 

%%%%%%%%%%%%%%%%%%%%%%%%%%%%%%%%%%%%%%%%%%%%%%%%%%%%%%%%%%%%%%%%%%%%
%\section*{1}
%%%%%%%%%%%%%%%%%%%%%%%%%%%%%%%%%%%%%%%%%%%%%%%%%%%%%%%%%%%%%%%%%%%%

\section{Introduction}

\vspace{0.5cm}

Vertex algebras naturally arise from highest weight representations of affine or Virasoro algebras. For a basic set up lets consider

\begin{equation}
D=t \frac{d}{dt}
\end{equation}

\noindent
acting as a differential operator on the ring $R=\mathbb{C}[t,t^{-1}]$. The Lie algebra $\mathfrak{g}$ of derivations of $R$ is generated by $t^mD, \ m \in \mathbb{Z}$. This algebra is graded, where $t^mD$ has weight $m$. We are interested to central extensions 

\begin{equation}
0 \to \mathbb{C}.c \to V \to \mathfrak{g} \to 0
\end{equation}

\noindent
namely the Virasoro algebra. It is also a graded Lie algebra. Certain highest weight representations of $V$ arise in conformal field theory. They have the property that for a suitable choice of $\widehat{D}$ lifting $D$ in the extension, the character 

\begin{equation}
Trace (q^{\widehat{D}})
\end{equation}

\noindent
is well defined and is the $q$-expansion of a modular form, \cite{BO}. The definition of Virasoro algebra can be extended to that of vertex algebra. We may think of the Virasoro algebra $V$ as a vector space acted on by the commuting operators $v_n$ and $V$ is generated by $1 \in V$ and the action of these operators. Then $V$ becomes a commutative ring such $1$ is the identity and the actions of all the operators are given by multiplication of elements of $V$. We define an operator 

\begin{equation}
\phi(t) : v \longmapsto \sum_iD^i v .t^i/i!, \qquad V \to End(V)[t,t^{-1}]
\end{equation}

\noindent
called vertex operator. The maps 

\begin{equation}
Trace (\phi(x)\phi(y)...)
\end{equation}

\noindent
are called correlation functions from their analogues in quantum field theory. A vertex operator algebra structure is expected to explain a conformal infinitesimal deformation of $V$. 

A vertex operator algebra is given by a $4$-tuple $(V,Y,1,\omega)$ where $V$ is a $\mathbb{Z}$-graded vector space with a linear map

\begin{equation}
Y(.,z): V \to End(V)[[z,z^{-1}]], \qquad Y(v,z)=\displaystyle{\sum_{n \in \mathbb{Z}}}v_n z^{-n-1}
\end{equation}

\noindent
where $1 \in V_0$ is called the vaccum vector with $Y(1,z)=id_V,\ v_{-1}1=v$. The vector $\omega \in V_2$ is a specific element called conformal element or Virasoro element such that

\begin{equation}
Y(\omega,z)=\sum_n L(n)z^{-n-1}
\end{equation}

\noindent
provides a set of Virasoro generators $L(n)$ such that $L(0)\mid_{V_n}=n.id_{V_n}$. The operator $L(-1)$ satisfies 

\begin{equation}
[L(-1),Y(v,z)]=\frac{d}{dz}Y(v,z)
\end{equation}

\noindent
and we also assume the Jacobi identity for the vertex operator $Y$, cf. \cite{M}. We will consider vertex algebras of CFT-type, that is $V_n=0, \ n<0$ and $V_0=\mathbb{C}.1$. If $V$ is generated by the subset $S \subset V$, then 

\begin{equation}
V=\text{span} \{v_{n_1}^1...v_{n_k}^k .\textbf{1}| v_i \in S\}
\end{equation} 

A unitary vertex operator algebra is one with a positive definite Hermitian form. This notion can also be defined for the modules over these algebras via an anti-involution. An anti-linear automorphism is a linear map such that

\begin{equation}
\phi:V \to V, \qquad \phi(1)=1, \ \phi(\omega)=\omega, \ \phi(u_n.v)=\phi(u)_n\phi(v)\qquad \forall u,v \in V 
\end{equation}

In contrast to Lie algebras the definition of invariant bilinear form on a vertex algebra is highly complicated. For instance, the contragrediant $(V',Y')$ of a vertex algebra module $(V,Y)$ is defined via the form

\begin{equation}
\langle Y'(v,x)w',w\rangle=\langle w', Y(e^{xL(1)}(-x^2)^{L(0)}v,x^{-1})w\rangle, \qquad v,w \in V,\ w' \in V'
\end{equation}

\noindent
The form should also explain the invariant bilinear form on our vertex algebra. A unitary vertex algebra is one which is equipped with a positive definite hermitian form. For a unitary vertex operator algebra the positive definite Hermitian form 

\begin{equation}
(.,.)_{\text{unitary}}:V \times V \to \mathbb{C}, \qquad \exists \lambda \in \mathbb{C}; \ (u,v)=\lambda (1,1)
\end{equation}

\noindent
is uniquely determined by its value at $(1,1)$.  This can be easily proved using the axioms of a vertex algebra and the invariance property in (11). 

Vertex algebras and vertex homomorphisms build up a tensor category. It simply means that we can tensor a finite number of vertex algebras 

\begin{equation}
\bigotimes_{i=1}^p(V_i,Y_{V_i},I_i,w_i), \qquad w=w_1 \otimes 1...\otimes 1 +...+1 \otimes ... 1\otimes w_p
\end{equation}

\noindent
The data of a vertex algebra should satisfy a locality axiom. Unfortunately the detailed description of this concept is out of the scope of this short note. We refer the reader to \cite{FB} for a complete explanation. It breifly means as follows. For any $A, B \in V$ the two formal power series in two variables obtained by composing $Y(A,z)$ and $Y(B,w)$ in two possible ways are equal, possiblly after multiplying them by a large power of $(z-w)$. It can be stated as

\begin{equation}
(z-w)^N[Y(A,z),Y(B,w)]=0, \qquad \text{some} \ N \in \mathbb{Z}_+
\end{equation}

\noindent
Regarding this we define a normally ordered product as

\begin{equation}
:A(z)B(w): \ =\sum_m \{\sum_{m<0}A_mB_nz^{-m-1}+\sum_{m \geq 0}A_mB_nz^{-m-1}\} w^{-n-1}
\end{equation}

\noindent
of vertex operators. The product can be inductively extended for more than two factors.

\begin{remark}\cite{B}
Vertex algebras when the operators $V(a,x)$ are holomorphic are commutative rings with derivations. The notation $V(u,z)v$ is a deformation of the one $u^z.v$. If we have a commutative algebra with a derivation $D$, then we can define 

\begin{equation}
V(a,x)b =\sum_{i \geq 0} (D^ia)bx^i/i!
\end{equation}

\noindent
Conversely if $V$ is a vertex algebra, we can define $ab=V(a,0)b$, and 

\begin{equation}
Da=\text{coefficient of}\ x^1 \ \text{in} \ V(a,x)b
\end{equation}

\noindent
Then in the new notation we put $a^x=\sum_ix^iD^ia/i!$, then 

\begin{equation}
a^xb=\sum_ix^iD^iab/i!
\end{equation}

Here $x$ is thought to be an element in the formal group $\hat{G}_a$. This formal group has its formal group ring, the algebra of polynomials $H=\mathbb{C}[D]$ and its coordinate ring is the ring of formal power series $\mathbb{C}[[x]]$. The tensor category of modules with a derivation is the same as the category of modules over the formal group ring $H$. So holomorphic vertex algebras are the same as commutative ring objects in this category. In the non-holomorphic case the expressions $a^xb^y...,\ (a,b \in V)$ are no longer holomorphic and can have singularities. In the new notation the identities of vertex algebra theory are easier to understand; for example

\begin{equation}
V(a,x)b=e^{xL_{-1}}(V(b,-x)a) \qquad   \Leftrightarrow \qquad a^xb=(b^{x^{-1}}a)^x
\end{equation}

\end{remark}

An interwining operator between 3 modules $(W_1,Y_1), \ (W_2,Y_2)$ and $(W_3,Y_3)$ is a linear map

\begin{equation}
I(.,z):W_2 \to Hom(W_3,W_1)\{z\}, \qquad u \mapsto I(u,z)=\sum_{n \in \mathbb{Q}} u_n z^{-n-1}
\end{equation}

\noindent
satisfying certain conditions of compatibility. The vertex operator  $Y_M(.,z)$ will become an interwining operator where $W_3=W_1=M$. One way to define these operators is 

\begin{equation}
I(w,z)v=e^{zL(-1)}Y_M(v,-z)w
\end{equation} 

The interwining operators are important tools in order to define a product structure in the category of vertex algebras. This is because the usual tensor product of two Lie algebra modules is not generally a module. This makes the definition of a product structure satisfying associativity in the tensor category of VOA's considerably complicated. Interwining operators also explain the base of the theory of conformal blocks. Using the $L(-1)$-property mentioned above and interwining operators one can extract a system of differential equations whose solution systems are these modules. These are $D$-modules which certain differential operators act on them, \cite{Y}. 

A polarized Hodge structure on a $\mathbb{Q}$-vector space $V$, is given by a representation

\begin{equation}
\phi: \mathbb{U}(\mathbb{R}) \to \text{Aut}(V_{\mathbb{R}}, \mathbb{Q}) , \qquad \mathbb{U}(\mathbb{R}) =\left( 
\begin{array}{cc}
a  &  -b\\
b &    a
\end{array} \right), \ a^2+b^2 =1
\end{equation}

\noindent
The group $G_{\mathbb{R}}=Aut(V_{\mathbb{R}},Q)$ is a real simple Lie group. The period domain $D$ associated to the Hodge structure $\phi$ is the moduli space of polarized Hodge structures on a fixed vector space $V$ with the same Hodge numbers. The group $G_{\mathbb{R}}$ acts transitively on the period domain $D$ by conjugation;

\begin{equation}
D=\{\phi:S^1 \to G_{\mathbb{R}} \ ; \ \phi=g^{-1}\phi_0 g \}
\end{equation}

\noindent
The isotropy group $H$ of a reference polarized Hodge structure $(V,Q,\phi)$ is a compact subgroup of $G_{\mathbb{R}}$, which contains a compact maximal torus $T$. The Lie algebra $\mathfrak{g}$ of the simple Lie group $G_{\mathbb{C}}$ is a $\mathbb{Q}$-linear subspace of $\text{End}(V)$, and the form $Q$ induces on $\mathfrak{g}$ a non-degenerate symmetric bilinear form $B: \mathfrak{g} \times \mathfrak{g} \to \mathbb{C}$ which upto scale is just the Cartan-Killing form $\text{tr}(\text{ad}(x)\text{ad}(y))$. For each point $\phi \in D$ 

\begin{equation}
\text{Ad}\phi:\mathbb{U}(\mathbb{R}) \to \text{Aut}(\mathfrak{g}_{\mathbb{R}},B)
\end{equation}

\noindent
is a Hodge structure of weight $0$ on $\mathfrak{g}$. This Hodge structure is polarized by $B$. Associated to each nilpotent transformation $N \in \mathfrak{g}$ one defines a limit mixed Hodge structure. The local system $\mathfrak{g} \to \Delta^*$ is then equipped with the monodromy $T=e^{\text{ad} N}$ and Hodge filtration defined with respect to the multi-valued basis of $\mathfrak{g}$ by $e^{\log(t)\frac{N}{2\pi i}}F^{\bullet}$, where $F^{\bullet}$ is the natural Hodge filtration on $\mathfrak{g}$ by (20). It gives a limit MHS $(\mathfrak{g},F^{\bullet},W(N)_{\bullet})$. The polarizing form gives perfect pairings 

\begin{equation}
B_k:Gr_k^{W(N)} \mathfrak{g} \times Gr_{-k}^{W(N)} \mathfrak{g} \to \mathbb{Q}, \qquad B_k(u,v)=B(v,N^kv)
\end{equation}

\noindent
defined via the hard Lefschetz isomorphism $N^k:Gr_{-k}^{W(N)}\mathfrak{g} \cong Gr_{k}^{W(N)}\mathfrak{g}$.

A family of projective manifolds defined via a proper smooth map 

\begin{equation}
f:X \to S \ , \qquad X_s=f^{-1}(s)
\end{equation}

\noindent
of quasi-projective varieties there associates a polarized variation of Hodge structures (VHS) 

\begin{equation}
(\mathcal{V}=R^kf_*\mathbb{C}, F^{\bullet})
\end{equation}

\noindent
If $\dim S=1$ this variation simply is understood as a topological deformation of the Hodge structure $V=H^k(X_s,\mathbb{C})$ over a punctured disc. The study of the asymptotic behavior of the VHS is an important issue in Hodge theory. We will always assume $V$ is equipped with the limit Hodge filtration. By the Riemann-Hilbert correspondence a local system of Hodge structures defines a $D$-module with flat connection on the base manifold $S$. This gives rise to a decreasing filtration denoted also by $F^{\bullet}=(F^i)$ on the vector bundle $\mathcal{V} \otimes_{\mathbb{Q}} \mathcal{O}_{S}$ by holomorphic sub-bundles, and a flat connection 

\begin{equation}
\nabla: \mathcal{V} \otimes_{\mathbb{Q}} \mathcal{O}_{S} \to \mathcal{V} \otimes_{\mathbb{Q}} \Omega_{S}^1
\end{equation}

\noindent
satisfying Griffiths transversality;

\begin{equation}
\nabla(F^i \mathcal{V}) \subset F^{i-1}\mathcal{V} \otimes \Omega_{S}^1
\end{equation}

These data are also equipped with a flat bilinear pairing 

\begin{equation}
P: \mathcal{V} \times \mathcal{V} \to \mathbb{Q}
\end{equation}
We compare the context of Hodge structures with that of vertex algebras of conformal field theory. Vertex algebras appear as the highest weight representations of infinite dimensional Lie algebras. A correspondence between Higgs bundles and opers already is known as non-abelian Hodge theorem due to C. Simpson. The Beilinson-Bernstein localization (correspondence) also compares the context of variation of Hodge structure with that of highest weight modules over flag manifolds of semisimple Lie groups. A more general analogue of the Bernstein correspondence over a local manifold can also be formulted in the context of geometric Langlands correspondence. We discuss a generalized version of Harish-Chandra modules called Wakimoto modules and a generalized Harish-Chandra homomorphism. We conclude with an explanation of geometric Langlands correspondence.

\vspace{0.3cm}

\textbf{Explanation on the text:} Section 1 is the introduction and we introduce the concept from the literature. 

In Section 2 we present main examples of vertex algebras we are dealing with, as affine Kac-Moody algebras and Virasoro algebras and Fock modules. We present Fock representations of Heisenberg algebra and Harish-Chandra pairs in this section. 

In Section 3 we give basic concepts related to variation of (mixed) Hodge structure. We explain the context of mixed Hodge modules and the non-abelian Hodge theorem of C. Simpson as equivalent notions. 

Section 4 contains the Beilinson-Bernstein localization functor which we successively develop over the $\mathfrak{g}$-opers in order to explain the geometric Langlands correspondence. We give a brief explanation of KZ-equations and the conformal blocks at the end.

\section{Vertex algebras}

In this section main examples of vertex algebras and their representations are presented along what we explained in the introduction. The references are \cite{B, DX, F1, F3, FB, FF, IK}. 

\begin{definition}
A vertex algebra consists of the following data;

\begin{itemize}
\item (space of states) A $\mathbb{Z}$-graded vector space 

\begin{equation}
V=\bigoplus_n V_n
\end{equation}

\item (vacuum vector) a vector $|0 \rangle \in V_0$
\item (translation operator) a linear operator $T:V \to V$ of degree one.
\item (vertex operators) a linear operation 

\begin{equation}
Y(.,z):V \to \text{End}(V)[[z^{\pm 1}]]
\end{equation}

\noindent
taking $A \in V_m$ to 

\begin{equation}
Y(A,z)=\sum_n A_{(n)} z^{-n-1}
\end{equation}

\noindent
of conformal dimension $m$, i.e $\deg A{(n)}=-n+m=1$.
\item (vacuum axiom) $Y(|0 \rangle )=Id_V$. Furthermore 

\begin{equation}
Y(A,z) |0 \rangle \in V[[z]] , \qquad \forall A \in V
\end{equation}

\item (translation axiom) For any $A \in V$,

\begin{equation}
[T,Y(A,z)]=\partial_z Y(A,z)
\end{equation}

\noindent
and $T|0 \rangle =0$.
\item (locality axiom) All fields are local with respect to each others.
\end{itemize}
\end{definition}

Vertex algebra structure naturally appear in known geometric concepts we may know. Lets begin with Lie algebra $H$ defined as central extension 

\begin{equation}
0 \to \mathbb{C}.1 \to H \to \mathbb{C}((t)) \to 0
\end{equation}

\noindent
It may also be regarded as the completion of the one dimensional central extension of the commutative Lie algebra of Laurent polynomials $\mathbb{C}[t,t^{-1}]$ having basis $b_n=t^n, \ n \in \mathbb{Z}$ and the central element $1$. Lets call the latter Lie algebra by $H'$. The universal enveloping algebra $U(H')$ is an associative algebra with generators $b_n$ and relations

\begin{equation}
b_nb_m-b_mb_n=n.\delta_{n,-m}1, \qquad b_n.1-1.b_n=0
\end{equation}

The left ideals $t^N \mathbb{C}[t]$ build up a system of open neighborhoods of $0$, and one can consider the completion of $U(H')$ with respect to this topology, denoted $\tilde{U}(H')$. The quotient 

\begin{equation}
\tilde{H}=\tilde{U}(H')/(1-1)
\end{equation}

\noindent 
is the well known Weyl algebra. Here the first $1$ is the central element and the second is the unit of $\tilde{U}(H')$. Let $\tilde{H}_+$ be the subalgebra of $\tilde{H}$ generated by $b_n, \ n \geq 0$ and define

\begin{equation}
V=Ind_{\tilde{H}_+}^{\tilde{H}} \mathbb{C}=\tilde{H}_-=\mathbb{C}[b_{-1},b_{-2},...]
\end{equation}

\noindent
The module $V$ is called the Fock representation of $\tilde{H}$. Now lets look at to the fields 

\begin{equation}
b(z)=\sum_n b_n z^{-n-1}
\end{equation}

\noindent
where $b_n$ is considered as an endomorphism of $V$. Since $\deg(b_n)=-n$, $b(z)$ is a field of conformal dimension one. Lets consider 

\begin{equation}
b(z)^2=\sum_n (\sum_{k+l=n}b_kb_l)z^{-n-2}
\end{equation}

\noindent
The relations (37) imply that the coefficient operators can be rearranged so that the annihilation operators ($b_n, \ n<0$) be in the right side of creation operators ($b_n, \ n \geq 0$) and for any $x \in V$ there are only a finite number of $b_kb_l$ whose action on $x$ is non zero. This makes the expression (41) well defined. There are standard ways in conformal field theory to remove infinite sums arising from repeatedly creating and annihilating the same state. In our case we define the normally ordered product of $b(z)$ with itself as 

\begin{equation}
:b(z)b(z)=\sum_n :b_kb_l:z^{-n-2}, \qquad :b_kb_l:=\begin{cases} b_lb_k \qquad l=-k, \ k \geq 0 \\ b_kb_l \qquad \text{otherwise} \end{cases}
\end{equation}

\noindent
With the normally ordered product we can proceed to define for instance 

\begin{equation}
Y(b_{-1}^2,z)=:b(z)^2:
\end{equation}

\noindent
etc .... The pattern explained above appears in many Lie algebra representations in finite or infinite dimensions. We will encounter several examples of this in the following.

\vspace{0.3cm}

\begin{itemize} 
\item[\textbf{(1)}] \textbf{Affine Kac-Moody algebras:} The first series of vertex operator algebras are affine Lie algebras. Let $\mathfrak{g}$  is a simple Lie algebra of finite dimensional over $\mathbb{C}$. Let $L\mathfrak{g}$ be the loop algebra $\mathfrak{g}((t))$. As a vector space an affine algebra is of the form $\widehat{\mathfrak{g}}:=L\mathfrak{g} \oplus \mathbb{C}.K$, with commutation relations $[K,.]=0$
and 

\begin{equation}
[A \otimes f(t),B \otimes g(t)]=[A,B] \otimes f(t)g(t)+(Res_{t=0}fdg(A,B)).K
\end{equation}

\noindent
where $(.,.)$ is an invariant bilinear form on $\mathfrak{g}$, normalized such that $(\theta,\theta)=2$ where $\theta$ is the highest root of $\mathfrak{g}$. 

A related concept is the vacuum representation of an affine algebra $\widehat{g}$. Let $k \in \mathbb{C}$, and suppose $\mathbb{C}_k$ be the 1-dimensional representation of $\mathfrak{g}[[t]] \oplus \mathbb{C}.k$. The vacuum representation of $\widehat{\mathfrak{g}}$ of level $k$ is 

\begin{equation}
V_k(\mathfrak{g})=Ind_{\mathfrak{g}[[t]] \oplus \mathbb{C}.K}^{\widehat{\mathfrak{g}}}\mathbb{C}_k
\end{equation}

\noindent
If $J^a$'s be a basis of $\mathfrak{g}$, then $J_n^a=J^a \otimes t^n$ and $K$ form a basis of $\widehat{\mathfrak{g}}$ and $V_k(\mathfrak{g})$ is generated by the monomials $J_{n_1}^{a_1}...J_{n_m}^{a_m}1_k$. We obtain a vertex algebra (module) $V_k(\mathfrak{g})$ with vertex operator

\begin{equation}
Y(J_{-1}^a .1_k,z)=J^a(z):=\sum_nJ_n^az^{-n-1}
\end{equation}

\noindent
The module $V_k(\mathfrak{g})$ has a unique maximal proper $\widehat{\mathfrak{g}}$-submodule $J(k)$ and 

\begin{equation}
L_{\mathfrak{g}}(k,0)=V_k(\mathfrak{g})/J(k)
\end{equation}

\noindent
becomes a simple vertex algebra. We denote by $L_{\mathfrak{g}}(k,\lambda)$ the corresponding highest weight module for $\widehat{\mathfrak{g}}$ associated to the highest weight $\lambda \in \mathfrak{h}^*$ of $\mathfrak{g}$, \cite{DX}, \cite{F1}.

\vspace{0.3cm}

\item[\textbf{(2)}] \textbf{Virasoro algebras:} The second type of vertex algebras we consider are Virasosro algebras denoted $\text{Vir}$. It is a central extension of the lie algebra $Der \mathbb{C}((t))$ generated by the operators $L_n=-t^{n+1}d/dt,\ n \in \mathbb{Z}$ by the 1-dimensional vector space $\mathbb{C}.C$ with relations, $[C,.]=0$ and 

\begin{equation}
[L_n.L_m]=(n-m)L_{n+m}+\frac{n^3-n}{12}\delta_{n,-m}C 
\end{equation}

\noindent
Assume $c, h \in \mathbb{C}$ and Let $\mathbb{C}_c$ be the one dimensional representation of $\text{Vir}$ defined by 

\begin{equation}
\begin{split}
L_n.1 &=0 \qquad \ n \geq 1\\
L_0.1 &=h.1\\
c.1 &=c.1
\end{split} 
\end{equation}

 Define 

\begin{equation}
V(c,h):=Ind_{Der\mathbb{C}[[t]] \oplus C}^{Vir} \mathbb{C}_c, \qquad c \in \mathbb{C}
\end{equation}

\noindent
Then $V(c,h)$ is generated by the monomials $L_{j_1}...L_{j_m}1_c, \ j_1 \leq j_2 \leq ... \leq -2$ and it is a vertex algebra module with vertex operator

\begin{equation}
Y(L_{-2}1_c,z)=T(z)=\sum_nL_n z^{-n-2}
\end{equation}

\noindent
The modules over the Virasoro algebra are classified according to the action of the operator $L(0)$ by $L(0).1=h.1$ and it becomes a highest weight module of the Virasoro algebra denoted $V(c,h)$. Its unique irreducible quotient is denoted by $L(c,h)$, \cite{DX}, \cite{IK}.

\vspace{0.3cm}

\item[\textbf{(3)}] \textbf{(Bosonic) Fock representations:}  Fock modules are defined using Heisenberg lie algebra of rank 1. 

\begin{equation}
H=\bigoplus_n \mathbb{C} a_n \oplus \mathbb{C}K, \qquad [K,.]=0,\ [a_m,a_n]=m\delta_{m+n,0}K
\end{equation}

\noindent
Let $\mathbb{C}_{\eta}$ be the 1-dimensional $H^{\geq}=(H^0 \oplus H^+)$-representation with 

\begin{equation}
a_n 1_{\eta}=\eta \delta_{n,0}1_{\eta}, \qquad K. 1_{\eta}=1_{\eta}
\end{equation}

\noindent
It has a $\mathbb{Z}$-gradation with 

\begin{equation}
\mathbb{C}_{\eta}^{n}=\begin{cases} \mathbb{C}_{\eta} \qquad n=0\\
0 \ \   \qquad n \ne 0 \end{cases}
\end{equation}

\noindent
The corresponding (bosonic) Fock module is defined by 

\begin{equation}
\mathcal{F}^{\eta}=Ind_{H^{\geq}}^H \mathbb{C}_{\eta}
\end{equation}

\noindent
The highest weight vector $1_{\eta}$ is denoted $|\eta \rangle$. We define a $\mathbb{Z}$-graded vertex algebra structure on $\mathcal{F}^0$ with vacuum vector $|0\rangle$, translation operator 

\begin{equation}
T|0\rangle=0,\ [T,a_n]=-na_{-n-1}
\end{equation}

\noindent
and vertex operator,

\begin{equation}
Y(a_{-1}|0\rangle, z)=\sum_n a_nz^{-n-1}
\end{equation}

Fock representations can be understood as the smallest representations of the Weyl algebra, \cite{IK}, \cite{F1}.

\vspace{0.3cm}

\item[\textbf{(4)}] \textbf{Harish-Chandra modules:} A pair $(\mathfrak{g},K)$ where $\mathfrak{g}$ is a Lie algebra, $K$ is a Lie group, such that $\mathfrak{k}=Lie(K)$ and an action 

\begin{equation}
Ad:K \to \mathfrak{g}
\end{equation}

\noindent
compatible with the adjoint action of $K$ on $\mathfrak{k}$ is called a Harish-Chandra pair. A $(\mathfrak{g},K)$-action on a scheme $X$ is a homomorphism 

\begin{equation}
\rho:\mathfrak{g} \to \Theta_X
\end{equation}

\noindent
together with an action of $K$ on $X$ such that 

\vspace{0.2cm}

\begin{itemize}
\item[(1)] the differential of the $K$-action is the restriction of the action of $\mathfrak{g}$ on $\mathfrak{k}$. 
\item[(2)] $\rho(Ad(k)(a))=k \rho(a)k^{-1}$. 
\end{itemize}

\vspace{0.2cm}

A (Harish-Chandra) $(\mathfrak{g},K)$-module is a vector space $V$ with the aformentioned compatible actions. One can consider the vector bundle 

\begin{equation}
\mathcal{V}=X \times_K V
\end{equation}

\noindent
on the scheme $X$, which gives a flat connection on the trivial vector bundle $X \times V$ on $X$. 

A Harish-Chandra module $M$ can also be defined over Virasoro algebras. Via this generalization the Lie algebra $\mathfrak{g}$ has a generalized triangle decomposition 

\begin{equation}
\mathfrak{g}=\bigoplus_{\alpha \in \mathfrak{h}^*} \mathfrak{g}_{\alpha}
\end{equation}

In this case the $\mathfrak{g}$-module $M$ is assumed to be $\mathfrak{h}$-diagonalizable 

\begin{equation}
M= \bigoplus_{\lambda \in \mathfrak{h}^*} M_{\lambda}, \qquad \dim M_{\lambda} < \infty
\end{equation}

\noindent
where each weight space is finite dimensional. Here $\mathfrak{h}$ is a Cartan subalgebra in the triangle decomposition $\mathfrak{g}=\mathfrak{g}_- \oplus \mathfrak{h} \oplus \mathfrak{g}_+$. It is known that $M$ must be a direct sum of highest weight module $Ind_{\mathfrak{g}^{\geq}}^{\mathfrak{g}} \mathbb{C}_{\Lambda}$ or a lowest weight module $Ind_{\mathfrak{g}^{\leq}}^{\mathfrak{g}} \mathbb{C}_{\Lambda}$ (Verma modules) or an intermediate series defined by $V_{a,b}, \ a,b \in \mathbb{C}$, 

\begin{equation}
V_{a,b}=\bigoplus_n \mathbb{C}v_n, \qquad L_s.v_n=(as+b-n)v_{n+s}, \qquad C.v_n=0
\end{equation}

\noindent
where $L_s$ are generators of the Virasoro algebra, see \cite{F1} and \cite{IK} for details.

\vspace{0.3cm}

\item[\textbf{(5)}] \textbf{Jantzen filtration and Shapovalov form:} Let $(\mathfrak{g},\mathfrak{h})$ be a lie algebra pair, $\mathfrak{h}$ a Cartan subalgebra, with the anti-involution $\sigma: \mathfrak{g} \to \mathfrak{g}$. The context of this item is applicable to any $Q$-graded Lie algebra, 

\begin{equation}
\mathfrak{g}=\bigoplus_{\beta \in Q} \mathfrak{g}_{\beta}, \qquad \dim \mathfrak{g}_{\beta} < \infty
\end{equation} 

\noindent
where $Q$ is an abelian group. In our case $Q=\mathfrak{h}^*$ is the root lattice. In partucular the definitions are applicable to usual finite dimensional Lie algebras as well as their affine or Virasoro algebras, which are infinite dimensional. Write

\begin{equation}
\mathfrak{g}=\mathfrak{g}_- \oplus \mathfrak{h} \oplus \mathfrak{g}_+
\end{equation} 

\noindent
The universal enveloping algebra of $\mathfrak{g}$ is by definition the quotient of the tensor algebra $T(\mathfrak{g})=\bigoplus_n \mathfrak{g}^{\otimes n}$ by the ideal generated by $x \otimes y -y \otimes x -[x,y]$. By the Poincar\'e-Birkhoff-Witt theorem we have the decomposition 

\begin{equation}
U(\mathfrak{g})=U(\mathfrak{h}) \oplus \{ \mathfrak{g}_-U(\mathfrak{g}) +U(\mathfrak{g})\mathfrak{g}_+ \}
\end{equation}

\noindent
Consider the projection 

\begin{equation}
\pi:U(\mathfrak{g}) \to U(\mathfrak{h}) \cong S(\mathfrak{h})
\end{equation}

\noindent 
with respect to this decomposition ($S(\mathfrak{h})$ is the symmetric algebra of $\mathfrak{h}$). The bilinear form 

\begin{equation}
F:U(\mathfrak{g}) \times U(\mathfrak{g}) \to S(\mathfrak{h}), \qquad F(x,y)=\pi (\sigma(x)y)
\end{equation}

\noindent
is called the Shapovalov form of $\mathfrak{g}$. It is a symmetric and  contravariant form, that is

\begin{equation}
F(zx,y)=F(x \sigma(y))
\end{equation} 

\noindent
The decomposition (55) implies similar decomposition on the universal enveloping algebra 

\begin{equation}
U(\mathfrak{g})=\bigoplus_{\beta} U(\mathfrak{g})_{\beta}, \qquad U(\mathfrak{g})_{\beta}=\{ x \in U(\mathfrak{g})|[h,x]=\beta(h)x,\ \forall h \in \mathfrak{h}\}
\end{equation}

\noindent
Then one can show that for $\beta_1 \ne \beta_2 \in Q$

\begin{equation}
F(x,y)=0, \qquad x \in U(\mathfrak{g})_{\beta_1}, \ y \in U(\mathfrak{g})_{\beta_2}
\end{equation}

For each $\beta \in Q$ one can choose a basis of $U(\mathfrak{g})_{-\beta}$ namely $X_j, \ j \in I$. The determinant 

\begin{equation}
D_{\beta}=\text{det}(F(X_i,X_j))_{i,j \in I} \in S(\mathfrak{h})
\end{equation}

\noindent
is called a Shapovalov determinant of $\mathfrak{g}$. 

The notion of contravariant bilinear form can be defined for any $\mathfrak{g}$-module $M$. A basic fact in this context is; any highest weight module $M$ has a unique contravariant bilinear form $\langle.,.\rangle:M \otimes M \to \mathbb{C}$ up to a constant ($\langle g.x,y\rangle=\langle x, \sigma(g).y\rangle, \ g \in U(\mathfrak{g})$). This fact can be easily checked for instance on Verma modules $M(\lambda)$ as the unique irreducible quotients of these modules. In case of Verma modules $M_{\lambda}$, the radical of such form is the maximal proper submodule $J(\lambda) \subset M(\lambda)$. 

The Jantzen filtration of a $\mathfrak{g}$-module is defined based on the Shapovalov form on $U(\mathfrak{g})$. Let $R=\mathbb{C}[t]$ and $\phi:R \to \mathbb{C}$ be the canonical map. Let $\tilde{M}$ be a free $R$-module or rank $r$ with a nondegenerate symmetric bilinear form 

\begin{equation}
(.,.)_{\tilde{M}} :\tilde{M} \times \tilde{M} \to R
\end{equation}

\noindent
Set $M=\phi \tilde{M}=M \otimes_R R/tR$. Then $M$ admits a symmetric bilinear form 

\begin{equation}
(\phi v_1,\phi v_2)=\phi(v_1,v_2)_{\tilde{M}}
\end{equation}

\noindent
For $m \in \mathbb{Z}_{\geq 0}$ set 

\begin{equation}
\tilde{M}(m)=\{ v \in \tilde{M}| (v, \tilde{M})_{\tilde{M}} \subset t^m R \} \stackrel{i_m}{\hookrightarrow} \tilde{M}
\end{equation}

\noindent
Then set $M(m)=Im \phi (i_m)$. It follows that

\begin{equation}
M=M(0) \supset M(1) \supset ...
\end{equation}

\noindent
defines a filtration of $\mathbb{C}$-vector spaces called Jantzen filtration. This filtration enjoys of the following properties (cf. \cite{IK}),

\vspace{0.2cm}

\begin{itemize}
\item $\bigcap_m M(m)=0$
\item $M(1)=rad(.,.)$
\item there exists a symmetric bilinear form $(.,.)_m$ on $M(m)$ such that $rad(.,.)_m=M(m+1)$.
\end{itemize}

The procedure of defining the Jantzen filtration appears both in the context of vertex algebras and variations of Hodge structure. The Jantzen filtration will correspond to the weight filtration in local systems of mixed Hodge structure. 

\vspace{0.3cm}

\item[\textbf{(6)}] \textbf{Unitary (conformal) vertex operator algebras:} Let $(V,Y,1,w)$ be a vertex algebra and $\phi:V \to V$ be an antilinear involution. $(V,\phi)$ is called unitary if there exists a positive definite Hermitian form 

\begin{equation}
(.,.):V \times V \to \mathbb{C}
\end{equation}

\noindent
such that for $ a,u,v \in V$ 

\begin{equation}
(Y(e^{zL(1)}(-z^{-2})^{L(0)}a,z^{-1})u,v)=(u,Y(\phi(a),z)v),
\end{equation}

\vspace{0.1cm}

\noindent
where 

\begin{equation}
\ Y(w,z)=\sum_n L(n)z^{-n-2}
\end{equation}

In a unitary vertex operator algebra the positive definite Hermitian form is uniquely determined by $(1,1)$ via the properties of vertex algebra mentioned before. 

\noindent
In the Virasoro case $V(c,h)$, there exists a unique Hermitian form with 

\begin{equation}
(1_{c,h},1_{c,h})=1,\qquad (L_nu,v)=(u,L_{-n}v)
\end{equation}

\noindent
It is known that $V(c,h)$ is unitary iff $c\geq 1,\ h \geq 0$ or $c=c_m,\ h=h_{r,s}^m$ where

\begin{equation}
c_m=1-\frac{6}{m(m+1)}, \qquad h_{r,s}^m=\frac{r(m+1)-sm)^2-1}{4m(m+1)}
\end{equation}

\vspace{0.1cm}

\noindent
In the affine case $V_{\mathfrak{g}}(k,\lambda)$ has a unique positive definite Hermitian form such that 

\begin{equation}
(1,1)=1, \qquad (xu,v)=-(u,\widehat{\omega}_0(x)v), \ x \in \widehat{\mathfrak{g}}, \ u,v \in L_{\mathfrak{g}}(k,\lambda)
\end{equation}

\noindent
where 

\begin{equation}
\widehat{\omega}_0:\widehat{\mathfrak{g}} \to \widehat{\mathfrak{g}}
\end{equation}

\noindent
is the Cartan involution. Then $V_{\mathfrak{g}}(k,\lambda)$ with $k \ne -h^{\vee}$ is a unitary vertex algebra iff $k \in \mathbb{Z}^+$ and $\lambda$ is a dominant integral weight satisfying $(\lambda,\theta)\leq k$. 

In the Fock module case we reduce to $M(1,\lambda)=U(\hat{\mathfrak{h}})/J_{\lambda}$. It is known that it is unitary iff $(\alpha,\lambda)\geq 0$, i.e $\lambda$ be a dominant weight, \cite{DX}. 

\vspace{0.3cm}

\item[\textbf{(7)}] \textbf{Conformal vertex algebras:} A vertex algebra $V=\oplus_n V_n$ of central charge $c$ is called conformal it it contains a vector $w \in V_2$ (called conformal vector) such that the corresponding vertex operator $Y(w,z)=\sum_n L_n z^{-n-2}$ satisfies 

\begin{equation}
L_{-1}=T,\qquad L_0|_{V_n}=n.Id, \qquad L_2 w=\frac{1}{2}c|0\rangle
\end{equation}

\noindent
It follows that there is a homomorphism 

\begin{equation}
Vir_c \to V, \qquad L_{-2}1_c \mapsto w
\end{equation}

\noindent
In the Kac-Moody case the conformal vector (called Sugawara conformal vector) is given by 

\begin{equation}
\frac{1}{2(k+h^{\vee})}\sum_a{(J_{-1}^a)}^2 1_k
\end{equation}

\noindent
where $J^a$ is an orthonormal basis of $\mathfrak{g}$. Thus a Kac-Moody algebra is conformal iff $k \ne -h^{\vee}$. In this case $\widehat{\mathfrak{g}}$ is a module over Virasoro algebra, \cite{F1}, \cite{F2}.

\end{itemize}

\section{Variation of Hodge structure}

The references for this section are \cite{DW, P, R, S, SA, SV}. A variation of Hodge structure over a complex manifold $S$ gives rise to a period map

\vspace{0.2cm}

\begin{center}
$\Phi:S \to \Gamma_{\mathbb{Z}} \backslash D$
\end{center}

\vspace{0.2cm}

\noindent
where $S$ is a smooth base manifold and $\Gamma$ is a discrete group. $D$ is the period domain and it is known that it is a hermitian symmetric complex manifold. There are naturally defined Hodge bundles $F^p$ of the Hodge structure on $V$, and also the endomorphism bundle associated to $\mathfrak{g}=End(V)$ on $D$. The corresponding local systems are  

\begin{equation}
\mathcal{V}:=\Gamma \backslash (D \times \mathcal{V}), \qquad \mathcal{G}:=\Gamma \backslash (D \times \mathfrak{g})
\end{equation}

\noindent
respectively. One way to explain the complex structure on $D$ is to embed it in its compact dual $\check{D}$, which is the set of all Hodge filtrations on $V$ with the same Hodge numbers satisfying the first Riemann-Hodge bilinear relation. $\check{D}$ is a homogeneous complex manifold. There are $G_{\mathbb{C}}$-homogeneous vector bundles 

\begin{equation}
F^p \to \check{D}
\end{equation}

\noindent
called Hodge bundles whose fiber at a given point $F^{\bullet}$ is $F^p$. Over $D \subset \check{D}$ we have $
V^{p,q}=F^p/F^{p+1}$, which are homogeneous vector bundles for the action of $G_{\mathbb{R}}$. They are Hemitian vector bundles with $G_{\mathbb{R}}$-invariant Hermitian metric given in each fiber by the polarization form. The space of functions on $D$ can be identified with the $\Gamma_{\mathbb{Z}}$-automorphic functions on $D$.

\vspace{0.3cm}
 
\begin{itemize}

\item[\textbf{(1)}] \textbf{Variation of mixed Hodge structure:} 
A polarized variation of mixed Hodge structure over the punctured disc $\Delta^*$ consists of the 5-tuple $( \mathcal{V}, F^{\bullet}, W_{\bullet}, \nabla, P)$ where 

\begin{itemize}

\item $\mathcal{V}$ is a local system of $\mathbb{Q}$-vector spaces on $\Delta^*$.
\item $W_{\bullet}$ is an increasing filtration on $\mathcal{V}$ by sub-local systems of $\mathbb{Q}$-vector spaces.
\item $F^{\bullet}=(F^i)$ is a decreasing filtration on the vector bundle $\mathcal{V} \otimes_{\mathbb{Q}} \mathcal{O}_{\Delta^*}$ by holomorphic sub-bundles.
\item $\nabla: \mathcal{V} \otimes_{\mathbb{Q}} \mathcal{O}_{\Delta^*} \to \mathcal{V} \otimes_{\mathbb{Q}} \Omega_{\Delta^*}^1$ is a flat connection satisfying Griffiths transversality;

\begin{equation}
\nabla(F^i) \subset F^{i-1} \otimes \Omega_{\Delta^*}^1
\end{equation}

\item $P: \mathcal{V} \times \mathcal{V} \to \mathbb{Q}$ is a flat pairing inducing a set of rational non-degenerate bilinear forms $P_k:Gr_W^k \mathcal{V} \otimes Gr_W^k \mathcal{V} \to \mathbb{Q}$ such that the triple 

\begin{equation}
(Gr_W^k \mathcal{V}, F^{\bullet}Gr_W^k ,P_k)
\end{equation}

\noindent
defines pure polarized variation of Hodge structure on $\Delta^*$. 

\end{itemize}

\noindent
which we briefly mention as $\mathcal{H}$.

Let $V$ a Hodge structure with an exhaustive decreasing Hodge filtration $F^p$. Regard a locally free sheaf $\xi(V,F)$ over $\mathbb{C}$ as the submodule of $V \otimes \mathbb{C}[t,t^{-1}]$ generated by $t^{-p}F^p$.  Given a real structure one can glue $\xi(V,F)$ and $\xi(V,\overline{F})$ using the involution $\ t \to (\bar{t})^{-1}$ to obtain a locally free sheaf $\xi(V,F,\overline{F})$ on $\mathbb{P}^1$ with the action of $\mathbb{C}^*$ and antilinear involution. This procedure may be explained as follows. A variation of polarized Hodge structure of weight $k$ provides a 4-tuple $(H,F, \nabla,P)$ where 

\begin{equation}
\nabla:H \to H \otimes z^{-1} \Omega_{\Delta^*}(\log 0)
\end{equation}

\noindent
is a flat connection and a $(-1)^k$-symmetric non-degenerate and flat pairing 

\begin{equation}
P:H \times j^* H \to \mathcal{O}_{\Delta^*}, \qquad j:z \mapsto -z
\end{equation}

\noindent
The bilinear form $P$ induces a non-degenerate symmetric pairing

\begin{equation}
z^{-k}P:H/zH \times H/zH \to \mathbb{C}
\end{equation}

The Hodge filtration can be explained as follows. Lets $V$ be the Kashiwara-Malgrange filtration on the mixed Hodge module associated to $H$. and suppose $(H,\nabla)$ is regular singular. Then $H \hookrightarrow V^{>-\infty}$. For $\alpha \in (0,1]$, define (cf. \cite{DW}),  

\begin{equation}
F^pH_{\lambda}:=z^{\beta+\frac{N}{2\pi i}}Gr_V^{-\beta}H
\end{equation}

We will identify the variation of polarized Hodge structures with their associated polarizable Hodge module via the Riemann-Hilbert correspondence. This correspondence has also been studied in a more systematic way by M. Saito \cite{SA}.

\vspace{0.3cm}

\item[\textbf{(2)}] \textbf{Polarizable Mixed Hodge modules:} 
Let $X$ be a complex algebraic variety and denote by $MHM(X)$, the abelian category of Mixed Hodge Modules on $X$. $MHM(X)$ is equipped with a forgetful functor 

\begin{equation}
\text{rat}:MHM(X) \to \text{Perv}(\mathbb{Q}_X)
\end{equation}

\noindent
which assigns the underlying perverse sheaf/$\mathbb{Q}$. Sometimes the above objects is understood as elements in $D^bMHM(X)$ and $D_c^b(\mathbb{Q}_X)$ respectively, and the same for the functor $\text{rat}$. When $X$ is smooth, then a mixed Hodge module on $X$ determines a 4-tuple $(M,F,K,W)$ where $M$ is a holonomic $D$-module with a \textit{good} filtration $F$ and, with rational structure 

\begin{equation}
\text{DR}(M) \cong \mathbb{C} \otimes K \in \text{Perv}(\mathbb{C}_X)
\end{equation}

\noindent
for a perverse sheaf $K$, and $W$ is a pair of weight filtrations on $M$ and $K$ compatible with $\text{rat}$ functor. $\text{DR}$ denotes the de Rham functor shifted by the $\dim(X)$. The de Rham functor is dual to the solution functor. If $X=pt$, Then, $MHM(pt)$ is exactly all the polarizable mixed Hodge structures.

A MHM always has a weight filtration $W$, and we say it is \textit{pure of weight $n$}, if $Gr_k^W=0$ for $k \ne n$. Normally, the filtration $W$ is involved with a nilpotent operator on $M$ or the underlying variation of a mixed Hodge structure. A mixed Hodge modules (def.) is obtained by successive extensions of pure one. If the support of a pure Hodge module as a sheaf is irreducible such that no sub or quotient module has smaller support, then we say the module has \textit{strict support}. Any pure Hodge module will have a unique decomposition into pure modules with different strict supports, known as Decomposition Theorem. A pure Hodge module is also called polarizable HM. $MH_Z(X,n)^p$ will denote the category of pure Hodge modules with strict support $Z$. 
An $M \in HM_Z(X,n)$ determines a polarizable variation of Hodge structure. The converse of this fact is also true, that variation of Hodge structures determine a MHM. Thus;

\begin{equation}
MH_Z(X,n)^p \simeq VHS_{gen}(Z,n-\dim Z)^p
\end{equation}

\noindent
The right side means polarizable variations of Hodge structure of weight $n-\dim Z$ defined on a non-empty smooth sub-variety of $Z$. Equation (97), explains a deep non-trivial fact about regular holonomic $D$-modules, their underlying perverse sheaves and their polarizations. It may also be interpreted as an analogue of Riemann-Hilbert correspondence between mixed Hodge modules and their underlying perverse sheaves, \cite{R}.

\vspace{0.3cm}

\item[(3)] \textbf{Higgs bundles and non-abelian Hodge theorem:} Suppose $X$ is smooth and projective over $\mathbb{C}$. A harmonic bundle on $X$ is a $C^{\infty}$-vector bundle $E$ with differential operators $\partial$ and $\overline{\partial}$ and algebraic operators $\theta$ and $\overline{\theta}$ such that the following holds: There exists a metric $h$ so that $\partial+\overline{\partial}$ is a unitary connection and $\theta + \overline{\theta}$ is self adjoint. And if 

\begin{equation}
\nabla=\partial+\overline{\partial} +\theta +\overline{\theta}, \qquad \nabla''=\overline{\partial}+\theta 
\end{equation}

\noindent
then $\nabla^2=\nabla''^2=0$. With these conditions $(E,D)$ is a vector bundle with flat connection, and $(E,\overline{\partial},\theta)$ is a Higgs bundle, i.e a holomorphic vector bundle with holomorphic section 

\begin{equation}
\theta \in H^0(End(E) \otimes \Omega_X^1), \qquad \theta \wedge \theta =0
\end{equation}

A Higgs bundle is stable (resp. semistable) if for any coherent subsheaf $F \subset E$ preserved by $\theta$ we have 

\begin{equation}
\deg(F)/\text{rank}(F) <\deg(E)/\text{rank}(E) \qquad (\text{resp}. \ \leq)
\end{equation}

There is a natural equivalence between the categories of harmonic bundles on $X$ and semisimple flat bundles (or representations of $\pi_1(X)$. There  is also a natural equivalence between the categories  of harmonic bundles and direct sum of stable Higgs bundles with vanishing Chern class. The resulting correspondence between representations and Higgs bundles can be extended to an equivalence between the category of all representations of $\pi_1(X)$ and all semistable Higgs bundles with vanishing Chern classes. This statement is refereed to as the non-abelian Hodge theorem. 

There is a natural $\mathbb{C}^*$ action on the category of semistable Higgs bundles with vanishing Chern classes, denoted 

\begin{equation}
t:(E,\theta) \mapsto (E,t\theta)
\end{equation}

\noindent
The semistable representations which are fixed by this action are exactly complex variations of Hodge structure. A representation $\varrho$ of $\pi_1(X)$ is called rigid if any nearby representation is conjugate to it. The correspondence described above is continuous on the moduli of semisimple representations. It follows that if a semisimple representation is called rigid it must be fixed by $\mathbb{C}^*$ and it comes from a complex variation of Hodge structure. In this case there is a $\mathbb{Q}$-variation of HS $V_{\mathbb{Q}}$ such that $\varrho$ is a direct factor of the monodromy representation of $V_{\mathbb{Q}} \otimes \overline{\mathbb{Q}}$ (the monodromy is a sum of conjugates of $\varrho$), \cite{S}. 

Let $M_{Dol}(G), \ M_{DR}(G),\ M_B(G)$ denote the moduli space of Higgs bundles of degree zero, local systems, and representations of $\pi_1(X)$ respectively. We will denote the smooth loci of these varieties by the superscript reg, as $M_{Dol}^{reg}$,..., and we usually omit the script reg future on. The non-abelian Hodge theorem gives a diffeomorphism 

\begin{equation}
\tau:M_{Dol}(G) \cong M_{DR}(G)
\end{equation}

\noindent
The Riemann-Hilbert correspondence between bundles with integrable connection and representations provides isomorphisms 

\begin{equation}
M_{DR}(G) \cong M_B(G)
\end{equation}

\noindent
see also \cite{S}, \cite{D}. A systematic study on the inter-relation of between the Higgs fields of Higgs bundles and the system of Hodge bundles in VHS can be found in \cite{P}. The corollary is a unipotent variation of mixed Hodge structure defines a Higgs field $\theta$ which is flat relative to $\nabla$ and $\partial+\bar{\partial}$. In this case the invariance under the $\mathbb{C}^*$-action in (101) explains the complex variation of Hodge filtration.

\end{itemize}

\section{Connection between Hodge structure and Vertex algebras}

The references of this section are \cite{Y, SV, K, FB, FF, F4, IK, M, F3, F2, F1, D, B, BD}. Let $\mathfrak{g}_{\mathbb{C}}=\mathfrak{g}_{\mathbb{R}} \otimes \mathbb{C}=Lie(G_{\mathbb{R}} \otimes \mathbb{C})$ be a complex semi-simple Lie algebra, $\mathfrak{h}=\mathfrak{t} \otimes \mathbb{C}$ a Cartan subalgebra, $\mathfrak{t}=Lie(T)$, and $K_{\mathbb{C}}$ a complex Lie group corresponding to the unique maximal compact subgroup $K \subset G_{\mathbb{R}}$. We denote $U(\mathfrak{g}_{\mathbb{C}})$ to be the universal enveloping algebra of $\mathfrak{g}_{\mathbb{C}}$. We assume the action of $K_{\mathbb{C}}$ will be locally finite, and its differential agrees with the corresponding subspace of $U(\mathfrak{g}_{\mathbb{C}})$. One may match these data with the case $D=G_{\mathbb{R}}/H$ is a general period domain sitting in the diagram

\begin{equation}
\begin{CD}
G_{\mathbb{R}}/T @>>> G_{\mathbb{C}}/B\\
@VVV      @VVV\\
D=G_{\mathbb{R}}/H @>>> G_{\mathbb{C}}/P=\check{D}
\end{CD}
\end{equation}

\vspace{0.2cm}

\noindent
with $T$ a maximal torus, $B$ a Borel subgroup, and horizontal arrows to be inclusions. Let $U(\mathfrak{h}_{\mathbb{C}})$ be the universal enveloping algebra of $\mathfrak{h}$. The Weyl group $W$ of $(\mathfrak{g}_{\mathbb{C}},\mathfrak{h}_{\mathbb{C}})$ acts on $U(\mathfrak{h}_{\mathbb{C}})$ and gives an isomorphism 

\begin{equation}
HC:Z(\mathfrak{g}_{\mathbb{C}}) \stackrel{\cong}{\longrightarrow} U(\mathfrak{h}_{\mathbb{C}})^W 
\end{equation}

\noindent
where the upper-index means the elements fixed by $W$. Using the isomorphism (105) one can assign to each positive root $\mu \in \mathfrak{h}_{\mathbb{C}}^*$ the homomorphism 

\begin{equation}
\chi_{\mu}: Z(\mathfrak{g}_{\mathbb{C}}) \to \mathbb{C}, \qquad z \mapsto HC(z)(\mu)
\end{equation}

\noindent
is called the infinitesimal character associated to $\mu$. A result of Harish-Chandra says that any character of $Z(\mathfrak{g}_{\mathbb{C}})$ is an infinitesimal character, and 

\begin{equation}
\chi_{\mu}=\chi_{\mu'} \qquad \Leftrightarrow \qquad \mu=w(\mu'), \ w \in W
\end{equation}

We will use the above set up in our construction of correspondence between Hodge bundles and vertex algebras. In fact our correspondence uses a generalization of (105). We will do this step by step as follows.

\vspace{0.3cm}

\begin{itemize}

\item[\textbf{(1)}] \textbf{Beilinson-Bernstein correspondence:} Let $G$ be a connected, complex reductive algebraic group $G$ defined over $\mathbb{R}$. Let $\mathfrak{b} =\mathfrak{h} \oplus \mathfrak{n} \subset \mathfrak{g}=Lie(G)$ be the Borel subalgebra, with unipotent radical 

\begin{equation}
\mathfrak{n}=\oplus_{\alpha \in \Phi^+}\mathfrak{g}^{-\alpha}
\end{equation}

\noindent
and Cartan algebra $\mathfrak{h}$. Any $\lambda$ in the weight lattice $\Lambda \subset \mathfrak{h}_{\mathbb{R}}^*$ lifts to an algebraic character $e^{\lambda}$ of the Borel subgroup $B \subset G$ corresponding to $\mathfrak{b}$. To this character there corresponds a unique $G$-equivariant line bundle $\mathcal{L} \to X=G/B$. The group $B$ acts as $e^{\lambda}$ on the fibers. Therefore $\Lambda$ can be regarded as the group of $G$-equivariant line bundles via $\lambda \mapsto \mathcal{L}_{\lambda}$. Let $\rho=\frac{1}{2}\sum_{\alpha \in \phi^+} \alpha \in \frac{1}{2}\Lambda$. Then 

\begin{equation}
\mathcal{L}_{-2\rho} \cong \bigwedge^nT^*X
\end{equation}

\noindent
If $G$ is simply connected then $\rho \in \Lambda$ and $\mathcal{L}_{-2\rho}$ has a well defined square root. Define

\begin{equation}
D_{\lambda}=\mathcal{O}(\mathcal{L}_{\lambda-\rho}) \otimes D_X \otimes \mathcal{O}(\mathcal{L}_{\rho-\lambda})
\end{equation}

\noindent
The Lie algebra $\mathfrak{g}$ acts by infinitesimal translation on sections of $\mathcal{L}_{\lambda-\rho}$. Thus $\mathfrak{g} \hookrightarrow \Gamma D_{\lambda}$ which induces 

\begin{equation}
U(\mathfrak{g}) \rightarrow \Gamma D_{\lambda}
\end{equation}

\noindent
where the center $Z(\mathfrak{g})$ of $U(\mathfrak{g})$ acts via the infinitesimal character $\chi_{\lambda}$. Therefore we get a homomorphism 

\begin{equation}
U_{\lambda}(\mathfrak{g})=U(\mathfrak{g})/\ker(\chi_{\lambda}) \to \Gamma D_{\lambda}
\end{equation}

\noindent
which is compatible with degree filtration. The Beilinson-Bernstein theorem asserts that the map (112) is an isomorphism. Let $Mod(U_{\lambda}(\mathfrak{g}))_{fg}$ be the category of finitely generated $U_{\lambda}(\mathfrak{g})$-modules or equivalently the category of finitely generated $U_{\lambda}(\mathfrak{g})$-modules on which $Z(\mathfrak{g})$ acts via the character $\chi_{\lambda}$. Also let $Mod(D_{\lambda})_{coh}$ refer to the category of coherent $D_{\lambda}$-modules. Define the two functors

\begin{equation}
\Delta:Mod(U_{\lambda}(\mathfrak{g}))_{fg} \to Mod(D_{\lambda})_{coh}, \qquad \Delta(M)=M \otimes_{U_{\lambda}(\mathfrak{g})} D_{\lambda}
\end{equation}

\begin{equation}
\Gamma:Mod(D_{\lambda})_{coh} \to Mod(U_{\lambda}(\mathfrak{g}))_{fg}, \qquad \Gamma(M)=H^0(X,M)
\end{equation}

\noindent
A theorem by Beilinson and Bernstein asserts that when $\lambda$ is a regular and integrally dominant weight the above functors define an equivalence of categories 

\begin{equation}
Mod(U_{\lambda}(\mathfrak{g}))_{fg} \cong Mod(D_{\lambda})_{coh}
\end{equation}

\noindent
In the same case for $\lambda$, a similar isomorphism will hold between the category of Harish-Chandra modules

\begin{equation}
HC(\mathfrak{g},K)_{\lambda} \cong HC(D_{\lambda},K)
\end{equation}

\noindent
where in the right hand side we mean the category of $D_{\lambda}$-modules with a compatible $K$-action. 

If we consider the polarization of $D$-modules as a Hermitian duality, 

\begin{equation}
P:M \times \overline{M} \to C^{\infty}(X_{\mathbb{R}}), \qquad \text{bilinear over} \ D_{\lambda} \times \overline{D}_{-\lambda}
\end{equation}

\noindent
In our settings, $M \mapsto \overline{M}$ defines a bijection $Mod(D_{X,\lambda})_{\text{rh}} \cong Mod(D_{\overline{X},-\lambda})_{\text{rh}}$. We sometimes forget the subscripts $X$ or $\overline{X}$ and simply write $D_{\lambda}$ and $D_{-\lambda}$. The pairing $P$ should be understood as a pairing in the form $(\sigma,\tau)=\int_X \langle \sigma, \bar{\tau} \rangle$ where $\langle.,.\rangle$ is flat hermitian pairing on the underlying vector bundles, and it is $\mathfrak{u}_{\mathbb{R}}$-invariant, where $\mathfrak{u}$ is the compact form of $\mathfrak{g}$ defined via a Cartan involution. In the correspondence (113-114) the flat bilinear form $P$ in (117) corresponds to the Shapovalov form and the weight filtration to the Jantzen filtration, \cite{SV}.

\vspace{0.3cm}

\item[\textbf{(2)}] \textbf{Localization Functor:} We formulate the Beilinson-Bernstein correspondence in a slightly more general language. Let $G$ be a connected simple Lie group over $\mathbb{C}$. Assume $LG=G((t))$ is the Lie group of $L\mathfrak{g}=\mathfrak{g}((t))$. Let $X$ be a smooth projective algebraic curve over $\mathbb{C}$ and $p \in X$ a point. Let $P$ be a principal $G$-bundle over $X$. Let 

\begin{equation}
\mathfrak{g}_P=P \times_G \mathfrak{g}
\end{equation}

\noindent
be the vector bundle associated to the adjoint representation of $G$. Let $\mathfrak{g}_{out}^P$ be the Lie algebra of sections of $\mathfrak{g}_P$ around $p \in X$, and let $G_{out}$ be the Lie group of $\mathfrak{g}_{out}$. There is a natural embedding 

\begin{equation}
\mathfrak{g}_{out}^P \hookrightarrow L\mathfrak{g}
\end{equation}

\noindent
which can be lifted to $\mathfrak{g}_{out}^P \rightarrow \widehat{\mathfrak{g}}$. Denote by $\mathcal{O}^0$ the category of $\widehat{\mathfrak{g}}$-modules where the Lie subalgebra $\mathfrak{g}_{in}=\mathfrak{g}[[t]]$ acts locally finite. The modular functor assigns to a module 

\begin{equation}
M \longmapsto  M/\mathfrak{g}_{out}^P M
\end{equation}

\noindent
The dual space of $M/\mathfrak{g}_{out}^P M$ is called the space of conformal blocks. The localization functor of Beilinson and Bernstein assigns to $M$ a $D$-module on the homogeneous space 

\begin{equation}
\mathcal{M}=LG/G_{out}
\end{equation}

\noindent
For any integer $k$ define a line bundle $\mathcal{L}^k$ on $\mathcal{M}$ together with a homomorphism from $\widehat{\mathfrak{g}}$ to the Lie algebra of infinitesimal symmetries of $\mathcal{L}^k$. This gives a homomorphism from $U_k(\widehat{\mathfrak{g}})$ to the algebra $D_k$ of global differential operators on $\mathcal{L}^k$. Thus for any $\widehat{\mathfrak{g}}$-module $M$ of level $k$ we can define a left $D_k$-module  

\begin{equation}
\Delta(M)=D_k \otimes_{U(\widehat{\mathfrak{g}})}M
\end{equation}

\noindent
The fiber of $\Delta(M)$ at $P$ is $ M/\mathfrak{g}_{out}^P M$, \cite{F2}. The map $\Delta$ appearing here is an analogue of (113) and generalize that over affine algebras. The adjoint functor is the global section functor. The same as previous section the Shapovalov form of $M$ correspond to the flat Hermitian pairing in the right side. We will employ the map 

\begin{equation}
\Delta:\widehat{\mathfrak{g}} \text{-Mod} \longrightarrow D_{\mathcal{M}} \text{-Mod}
\end{equation}

\noindent
which is a generalization of its analogue (113), in the procedure of geometric Langlands correspondence in Session (6).  

\vspace{0.3cm}

\item[\textbf{(3)}] \textbf{$\mathfrak{g}$-opers on the punctured disc:} We first define $GL(n)$-opers, that is a pair $(F_{\bullet} \subset E,D)$, where $E$ is a vector bundle of rank $n$ equipped with 

\vspace{0.2cm}

\begin{itemize}
\item a complete flag $(0) \subset F_1 \subset ... \subset F_n =E, \ rank(F_i)=i$, \\[0.05cm]
\item a holomorphic connection $D:E \to E \otimes K$ necessarily flat such that, \\[0.05cm]
\item Griffiths transversality; $D:F_i \to F_{i+1} \otimes K$, \\[0.05cm]
\item non-degeneracy (strictness) $Gr_i D:Gr_i^F E \cong Gr_{i+1}^F E \otimes K$.
\end{itemize}

\noindent
An $SL(n)$-oper is a $GL(n)$-oper such that $\det(E)=\mathcal{O}$ and $D$ induces trivial connection $d$ on $\det E$. 

Now assume $G$ is a simple complex group, with $B \subset G$ a fixed Borel subgroup, $N=[B,B]$ its unipotent radical, $H=B/N$ the Cartan subgroup, and $Z=Z_G$ the center. The corresponding Lie algebras are $\mathfrak{n} \subset \mathfrak{b} \subset \mathfrak{g}, \mathfrak{t}=\mathfrak{b}/\mathfrak{n}$. If $P_B$ is a holomorphic principal $B$-bundle, $P_G$ the induced $G$-bundle and $Conn(P_G)$ the sheaf of connections. Define the projection 

\begin{equation}
c:Conn(P_G) \to \mathfrak{g}/\mathfrak{b}_P \otimes K
\end{equation}

\noindent
by requiring $c^{-1}(0)=Conn(P_B)$ and $c(D+\nu)=c(D)+[\nu]$ where $\nu$ is a section of $\mathfrak{g}_P \otimes K$ and $[\nu]$ is its image in $\mathfrak{g}/\mathfrak{b}_P \otimes K$. We will consider the class 

\begin{equation}
c(D) \in H^0(\mathfrak{g}/\mathfrak{b}_P \otimes K)
\end{equation}

\noindent
by taking locally any flat $B$-connection $D_B$, and then glue the local sections $[D-D_B]$. Then, a $G$-oper on $X$ is a pair $(P_B,D), \ D \in H^0(Conn(P_G))$ such that

\begin{itemize}
\item $c(D) \in H^0((\mathfrak{g}_{-1})_P \otimes K) \subset  H^0(\mathfrak{g}/\mathfrak{b}_P \otimes K)$\\[0.05cm]
\item For any simple negative root $\alpha$ the component $c(D)_{\alpha} \in H^0(\mathfrak{g}/\mathfrak{b}_P \otimes K)$ is nowhere vanishing. 
\end{itemize}

\noindent
The meaning of the conditions is that the connection $D$ preserves the flag corresponded to the Borel subgroup $B$ via the Griffiths transversality. By definition a $\mathfrak{g}$-oper is a $G$-oper where $G$ is the group of inner automorphisms of $\mathfrak{g}$. 

For $GL(n)$ the oper condition implies that if $E_U \cong \mathcal{O}^n$ is a trivialization compatible with the flag on an open chart $U$, then one can write the flat connection 

\begin{equation}
d+\left( 
\begin{array}{ccccc}
*  &  * & ... & * & *\\
\times &  *  & ... & * & * \\
0 & \times & ... & * & *\\
... & ... & ... & ... & *\\
0 & 0 & ... & \times & *
\end{array} \right)dt
\end{equation}

\noindent
Then a $B$-Guage equivalence class of a $GL_n$-oper has a unique representative of the form 

\begin{equation}
d+\left( 
\begin{array}{ccccc}
a_1  &  a_2 & ... & ... & a_n\\
1 &  0  & ... & 0 & 0 \\
0 & 1 & ... & 0 & 0\\
... & ... & ... & ... & 0\\
0 & 0 & ... & 1 & 0
\end{array} \right)dt
\end{equation}

\noindent
When we discuss about $\mathfrak{g}$-opers we automatically are considering their $B$-guage equivalence classes, \cite{D}. Lets restrict to opers on the punctured disc. By definition the space of $\mathfrak{g}$-opers on the punctured disc $D^{\times}$ is 

\begin{equation}
\mathcal{O}\mathfrak{p}_{\mathfrak{g}}(D^{\times})=\{\sum_i \psi_i X_{-\alpha_i}+v \big| 0 \ne \psi_i \in \mathbb{C}((t)), \ v \in \mathfrak{b}((t))\}/B((t))
\end{equation} 

\noindent
where $\alpha_i$ are the set of positive simple roots of $\mathfrak{g}$ with respect to $B$. The action of $B((t))$ is via the guege transformation by

\begin{equation}
g. D=Ad(g)D-(\partial_tg)g^{-1}
\end{equation} 

An oper on the punctured disc is called nilpotent if its connection has regular singularity at the origin with unipotent monodromy. $\mathfrak{g}$-opers as guage equivalence classes of flat connections can be compared with mixed Hodge modules. 

\vspace{0.3cm}

\item[\textbf{(4)}] \textbf{$\widehat{\mathfrak{g}}$-modules associated to opers:} Let $\mathfrak{g}$ be a finite dimensional semisimple Lie algebra and $\widehat{\mathfrak{g}}$ its affine Lie algebra. Define 

\begin{equation}
\tilde{U}(\widehat{\mathfrak{g}}): = \lim_{\leftarrow}U(\widehat{\mathfrak{g}})/U(\widehat{\mathfrak{g}})(\mathfrak{g} \otimes t^n\mathbb{C}[t])
\end{equation}

\noindent
It follows from a formula of Kac-Kazhdan for the determinant of the Shapovalov form that the module $V_k(\widehat{\mathfrak{g}})$ contains null vectors other than the highest weight vector $v_k$ only if $k=-h^{\vee}$ (the critical level). The space of null vectors $\hat{\mathfrak{z}}(\mathfrak{g})$ of $V:=V_{-h^{\vee}}$ is isomorphic to $End_{\widehat{\mathfrak{g}}}(V)$. To each vector $v \in V$ we can associate a power series 

\begin{equation}
v \longmapsto Y(v,z)=\sum_m v_mz^m
\end{equation}

\noindent
as in (27). The coefficient of these power series are elements of $\tilde{U}_{-h^{\vee}}(\widehat{\mathfrak{g}})=\tilde{U}(\widehat{\mathfrak{g}})/(K+h^{\vee})$. They span a Lie subalgebra $\tilde{U}_{-h^{\vee}}(\widehat{\mathfrak{g}})_{loc}$. For example 

\begin{equation}
A \in \mathfrak{g} \longmapsto Y((A \otimes t^{-1})v,z)=A(z)= \sum_n (A \otimes t^n)z^{-n-1}
\end{equation}

\noindent
This shows $\widehat{\mathfrak{g}} \subset \tilde{U}_{-h^{\vee}}(\widehat{\mathfrak{g}})_{loc}$. Let ${Z}(\widehat{\mathfrak{g}})$ be the center of $\tilde{U}_{-h^{\vee}}(\widehat{\mathfrak{g}})_{loc}$. One can show that 

\begin{equation}
x \in \mathfrak{z}(\widehat{\mathfrak{g}})  \qquad \Leftrightarrow \qquad Y(x,z) \in Z(\widehat{\mathfrak{g}})
\end{equation}

\noindent
and all the elements of $Z(\widehat{\mathfrak{g}})$ can be obtained in this form. A basic example of this is the Cassimir element 

\begin{equation}
S=\frac{1}{2} \sum_n (J_a \otimes t^{-1})^2 \ \in \mathfrak{z}(\widehat{\mathfrak{g}})
\end{equation}

\noindent
where $J_a \ , a=1,..., \dim \mathfrak{g}$ is an orthonormal basis of $\mathfrak{g}$ with respect to the invariant bilinear form. The coefficient $S_n$ of the power series

\begin{equation}
Y(S,z)=\sum_n S_n z^{-n-2} = \frac{1}{2}\sum_a :J_a(z)^2:
\end{equation}

\noindent
are called the Sugawara operators and lie in $Z(\widehat{\mathfrak{g}})$. Let 

\begin{equation}
M_{\chi,k}=U_k(\widehat{\mathfrak{g}}) \otimes _{U(\tilde{\mathfrak{b}}_+)} \mathbb{C}_{\chi}, \qquad \tilde{\mathfrak{b}}_+ =(\mathfrak{b}_+ \otimes 1) \oplus (\mathfrak{g} \otimes t \mathbb{C}[[t]])
\end{equation}

\noindent
be a Verma module over $\widehat{\mathfrak{g}}$ and $\chi \in \mathfrak{h}^*$ with highest weight vector $v_{\chi,k}$. Denote by $\mathfrak{g}^L$ the Langlands dual of $\mathfrak{g}$ obtained by exchanging the role of roots and co-roots. By classical theorems the center $Z(\widehat{\mathfrak{g}})$ is isomorphic to $W(\mathfrak{g}^L)$; the space of local functionals on $\mathcal{O}\mathfrak{p}(\mathfrak{g}^L)$. Let $\rho \in \mathcal{O}\mathfrak{p}(\mathfrak{g}^L)$ be  $\mathfrak{g}^L$-oper on the punctured disc. Then by what was said $\rho$ defines a central character $\tilde{\rho}:Z(\widehat{\mathfrak{g}}) \to \mathbb{C}$. Then  we associate the $\widehat{\mathfrak{g}}$-module $M_{\chi}^{\rho}=M_{\chi,-h^{\vee}}/\ker \tilde{\rho}$ to the ${\mathfrak{g}}^L$-oper $\rho$, \cite{F2}. Therefore we obtain a map

\begin{equation}
\mathcal{O}\mathfrak{p}(\mathfrak{g}^L) \longrightarrow Mod(\widehat{\mathfrak{g}}), \qquad \rho \longmapsto M_{\chi}^{\rho}=M_{\chi,-h^{\vee}}/\ker \tilde{\rho}
\end{equation}

\noindent
This is part of the correspondence between the connections and $\widehat{\mathfrak{g}}$-modules. The assignment or interpretations of $\mathfrak{g}^L$-opers as characters on the center of the dual affine algebra $Z(\widehat{\mathfrak{g}})$ is quite fundamental in this context. 

\vspace{0.3cm}

\item[(5)] \textbf{Wakimoto modules:} Our explanation of Wakimoto modules is quite brief, however we need this step in order to complete our correspondence. Assume we are given a linear function $\chi:\mathfrak{h}((t)) \to \mathbb{C}$. We extend it trivially to $\mathfrak{n}((t))$ and obtain a linear function on $\mathfrak{b}_-((t))$ also denote it by $\chi$. Now instead of considering $Ind_{\mathfrak{b}_{-}((t))}^{\mathfrak{g}((t))}\mathbb{C}_{\chi}$ we will consider the semi-infinite induction \cite{FF}, \cite{F4}. The resulting module is a module over the central extension of $\mathfrak{g}((t))$ i.e over $\widehat{\mathfrak{g}}$ of critical level where the vacuum is killed by $t \mathfrak{g}[t] \oplus \mathfrak{n}_-$. The parameters of the module will no longer behave as functionals on $\mathfrak{h}((t))$, but as connections on the $H^L$-bundle $ \Omega^{i}$ (sheaf of $i$-forms). They are precisely elements of the space $Conn_{\mathfrak{g}^L}(D^{\times})$ ($D^{\times}$ is the punctured disc). We obtain a family of smooth representations of $\tilde{U}_{\kappa_c}\widehat{\mathfrak{g}}$ parametrized by $Conn_{\mathfrak{g}^L}(D^{\times})$. For $\chi \in Conn_{\mathfrak{g}^L}(D^{\times})$ denote the corresponding Wakimoto module by $W_{\chi}$. The center $Z_{\kappa_c}\widehat{\mathfrak{g}}$ acts on $W_{\chi}$ according to a character. The corresponding point in $Spec (Z_{\kappa_c}\widehat{\mathfrak{g}})=\mathcal{O}p_{\mathfrak{g}^L}(\Delta^*)$ is denoted by $\mu(\chi)$. We obtain a map 

\begin{equation}
\mu:Conn_{\mathfrak{g}^L}(D^{\times}) \to \mathcal{O}p_{\mathfrak{g}^L}(D^{\times})
\end{equation}

\noindent
called the Miura transformation. Regarding the context of opers in Section (3) the Miura transformation should be understood as

\begin{equation}
\partial_t-\left( 
\begin{array}{ccccc}
0  &  q_1(t) & ... &  & q_{n-1}(t)\\
1 &  0  & ... & 0 & 0 \\
0 & 1 & ... & 0 & 0\\
... & ... & ... & ... & 0\\
0 & 0 & ... & 1 & 0
\end{array} \right)dt \longmapsto 
\partial_t-\left( 
\begin{array}{ccccc}
\chi_1(t)  &  0 & ... & ... & 0\\
1 &  \chi_2(t)  & ... & 0 & 0 \\
0 & 1 & ... & 0 & 0\\
... & ... & ... & ... & 0\\
0 & 0 & ... & 1 & \chi_n(t)
\end{array} \right)dt
\end{equation} 

\noindent
This amounts to the following splitting of the differential operator

\begin{equation}
\partial_t^n-q_1(t)\partial_t^{n-1}-...-q_{n-1}(t)=(\partial_t-\chi_1(t))...(\partial_t-\chi_n(t))
\end{equation}

\noindent
The fiber of the Miura transformation over a nilpotent oper is the variety of all Borel subalgebras $\mathfrak{g}^L$ containing that oper. It is called Springer fiber over the nilpotent oper. For example the Springer fiber over $0$ is the flag variety of $\mathfrak{g}^L$, \cite{F3}. The composition of (137) and (138) gives a map

\begin{equation}
Conn_{\mathfrak{g}^L}(D^{\times}) \longrightarrow \mathcal{O}\mathfrak{p}(\mathfrak{g}^L) \longrightarrow Mod(\widehat{\mathfrak{g}}) 
\end{equation}

\noindent
which is the way we correspond $D$-modules of flat connections to $\widehat{\mathfrak{g}}$-modules or vertex algebras. The map in (141) is an analogue of Harish-Chandra homomorphism 

\begin{equation}
HC:Z(\mathfrak{g}) \to \mathbb{C}[\mathfrak{h}^*]^W
\end{equation}

\noindent
A variation of MHS on the punctured disc can be interpreted as a flat connection that is regarded as an element of $Conn_{\mathfrak{g}^L}(D^{\times})$ via the action of $\mathfrak{g}^L$ coming from the internal symmetries of the Hodge structure. In many cases the Lie algebra action is paired with a compatible action of the Lie group $K$ to build up a Harish-Chandra pair $(\mathfrak{g},K)$.  

\vspace{0.3cm}

\item[\textbf{(6)}] \textbf{Geometric Langlands correspondence:} A more solid way to reformulate the localization theorem of Beilinson-Bernstein is via the Geometric Langlands correspondence. We state this briefly here to provide the idea. Our reference is \cite{F2}. Let $X$ be a smooth projective curve defined over a finite field $\mathbb{F}_q$ and $G$ a split connected simple algebraic group defined also over $\mathbb{F}_q$. For a closed point $x \in X$ let $\mathcal{O}_x$ be the completion of the local ring of $x$; i.e. $\mathcal{O}_x={\mathbb{F}_q}_x[[t]]$ and $K_x$ its field of fractions, where $q_x=q^{\deg x}$.

\textit{There is a correspondence between the conjugacy classes in $G^L(\overline{\mathbb{Q}_l})$ and $G_x$-modules which is called the local Langlands correspondence.}

\textit{The global Langlands correspondence states that; an irreducible unramified representation $\otimes_x' \pi_x$ is automorphic iff there exists a continuous homomorphism} 

\begin{equation}
\sigma:\pi_1(X) \to G^L(\overline{\mathbb{Q}_l})
\end{equation}

\noindent
\textit{such that each $\pi_x$ corresponds to the conjugacy class $\sigma(Fr_x)$ in the case of local Langlands correspondence.}

Now lets $X$ is an algebraic curve defined over $\mathbb{C}$. Let $G_{in}=G[[t]]$ and $\mathcal{O}_{crit}^0$ be the category of unramified $\widehat{\mathfrak{g}}$-modules of critical level, which consist of modules on which the action of $\mathfrak{g}_{in}=\mathfrak{g}[[t]]$ is locally finite and which contain $\mathfrak{g}_{in}$-invariant vector. On such modules the action of $\mathfrak{g}_{in}$ can be integrated to an action of the Lie group $G_{in}$. The analogue of a conjugacy class in the group $G^L$ is a regular $\mathfrak{g}^L$-oper on the formal disc. The analogue of the local Langlands correspondence is the following: 

\textit{Each regular $\mathfrak{g}^L$-oper $\rho_x$ on the formal disc defines an irreducible $\widehat{\mathfrak{g}}$-module of critical level.} 

Suppose that we are given a $\mathfrak{g}^L$-oper $\rho_x$ for each point $x \in X$. Let $V^{\rho_x}$ be the corresponding $\widehat{\mathfrak{g}}$-module of critical level defined at the end of Session (4). Set $\mathfrak{g}(\mathbb{A})=\prod_x'\mathfrak{g}((t_x))$ and $\widehat{\mathfrak{g}}(\mathbb{A})$ be its one dimensional central extension. The product $\otimes_x'V^{\rho_x}$ is naturally a $\widehat{\mathfrak{g}}(\mathbb{A})$-module. In this case one can assign a twisted $D$-module to the $\widehat{\mathfrak{g}}(\mathbb{A})$-module $\otimes_xM_x$ by the localization functor defined before;

\begin{equation}
loc: \otimes_x'V^{\rho_x} \longmapsto \Delta(\otimes_x M_x)
\end{equation}

\noindent
The action of $\mathfrak{g}_{in}$ on $M_x$ can be integrated to an action of $G_{in}$. Therefore the $D$-module is $K$-equivariant and descends to a $D$-module on $\mathcal{M}$ in (79), denoted by $\Delta(\otimes_x M_x)$. Lets specialize to $k=-h^{\vee}$. The $\widehat{\mathfrak{g}}(\mathbb{A})$-module $\otimes_x V^{\rho_x}$ is called weakly automorphic if ${\Delta}(\otimes_x V^{\rho_x}) \ne 0$. The weak version of the global Langlands correspondence over $\mathbb{C}$ states as follows: 

\textit{The $\widehat{\mathfrak{g}}(\mathbb{A})$-module $\otimes_x V^{\rho_x}$ is weakly automorphic iff there exists a globally defined regular $\mathfrak{g}^L$-oper $\rho$ on $X$ such that for each $x$, the $\rho_x$ is the restriction of $\rho$ to a small disc around $x$.}

\vspace{0.3cm}

\item[\textbf{(5)}] \textbf{Knizhnik-Zamolodchikov (KZ)-equations:} We keep this section as a well known example in conformal field theory. Assume $\mathfrak{g}$ is a finite dimensional semisimple Lie algebra with invariant bilinear form $\kappa$. Let $\widehat{\mathfrak{g}}_k$ be the affine lie algebra with level $k$ and dual Coxeter number $h^{\vee}$. The null vector of a $\widehat{\mathfrak{g}}_k$-module defines differential equations

\begin{equation}
(k+2)\frac{\partial}{\partial z_i} \Psi=\sum_{i \ne j} \frac{\Omega_{ij}}{z_i-z_j}\Psi
\end{equation}

\noindent
called KZ-equations, where $\Omega_{ij}=\sum_aJ^aJ_a$ are matrices. $J^a$ and $J_a$ are dual basis with respect to the invariant bilinear form $\kappa$ on $\mathfrak{g}$. In case that $\mathfrak{g} $ is semisimple, its local system corresponds to representations of the braid group

\begin{equation}
\theta:B_n \to V_1 \otimes ...\otimes V_n
\end{equation}

\noindent
as the holonomy of the Hamiltonian system (145) (by Riemann-Hilbert correspondence or non-abelian-Hodge theorem). A more concrete way to write this equation is 

\begin{equation}
dw=\sum_{1 \leq i < j \leq n} \frac{dz_i-dz_j}{z_i-z_j}A_{ij}w
\end{equation}

\noindent
which can be written as an equation of a flat connection $\nabla^{KZ}=d-\Gamma$ where $\Gamma=\sum (dz_i-dz_j)A_{ij}/(z_i-z_j) $. In the simplest case of 3-correlations over $\mathbb{P}^1 \setminus 0,1,\infty$ the equation (147) can be reduced to a one variable equation

\begin{equation}
\frac{d\phi}{dx}=(\frac{A}{x}-\frac{B}{1-x})\phi, \qquad \phi \in W\{x\}[\log x]
\end{equation}

\noindent
after suitable change of variables, where $A, \ B \in End(W)$ for a $\mathfrak{g}$-module $W$ are diagonal matrices. The system of solutions 
for (148) can be described in two ways which we explain as follows

\begin{itemize}

\item Suppose $\{\lambda\} $ is the set of eigenvalues of $A$ such that all eigenvalues of $A$ are contained in $\cup_{\lambda} \ \lambda +\mathbb{N}$. Thus for each $\lambda $ there is a set $\{\lambda +N_j^{\lambda}\}_{j=0}^{J_{\lambda}}$ such that $0=N_0^{\lambda}  <...< N_{J_{\lambda}}^{\lambda}$. Denote by $\pi_{\lambda}^A$ the projection onto the $\lambda$-eigenspace. Then a basic but perhaps not short calculation shows that for any $w \in W$ there is a unique solution to (148) of the form

\begin{equation}
\phi_w^A(x)=\sum_{\lambda} \sum_j \sum_{i \geq N_j^{\lambda}}w_{i,j}^{(\lambda)}x^{\lambda+i}(\log x)^j
\end{equation}

\noindent
where $w_{0,0}^{(\lambda)}=\pi_{\lambda}^A(w)$ and for each $j >0$ we have $\pi_{\lambda_{N_j^{\lambda},0}}^A(w_{\lambda_{N_j^{\lambda},0}}^{(\lambda}))=\pi_{\lambda_{N_j^{\lambda},0}}^A(w)$.

\item With similar set up but this time looking at the $B$-eigenvalues one shows that; for any $w \in W$ there is a unique solution of 

\begin{equation}
\frac{d\phi}{dy}=(\frac{B}{y}-\frac{A}{1-y})\phi, \qquad \phi \in W\{y\}[\log y]
\end{equation}

\noindent
of the form 

\begin{equation}
\phi_w^B(y)=\sum_{\mu} \sum_k \sum_{i \geq M_k^{\mu}}w_{i,k}^{(\mu)}y^{\mu+i}(\log y)^k
\end{equation}

\end{itemize} 

In the first case the map 

\begin{equation}
\phi_A:w \mapsto \phi_w^A(z) 
\end{equation}

\noindent
defines a $\mathfrak{g}$-isomorphism between $W$ and the solution system of the KZ-equation. In the second case the map

\begin{equation}
\phi_B:w \mapsto \phi_w^B(1-z) 
\end{equation}

\noindent
The $\mathfrak{g}$-automorphism $\Phi_{KZ}=\phi_B^{-1} \phi_A$ of $W$ is called the Drinfeld associator of $W$, which can be interpreted by the interwining operators explaining the tensor structure on the solutions of KZ-equation, \cite{M}.

\vspace{0.3cm}

\item[\textbf{(6)}] \textbf{Conformal blocks:} We conclude with this section as a brief from \cite{FB} to be continued in the next papers. As was mentioned a vertex algebra is called conformal if among the vertex operators, there is the generating function of the basis element of the Virasoro algebra. Such algebras or their modules become a module over Virasoro algebra. To a conformal vertex algebra $V$ one can assign a vector bundle $\mathcal{V}$ over an algebraic curve as a base manifold $X$. A vertex operator in this set up may be interpreted as a section of the dual bundle $\mathcal{V}^*$ on the punctured disc $D_x^{\times}$ with values in $End \mathcal({V}_x)$, which may be written as

\begin{equation}
A \otimes z^ndz \longmapsto Res_{z=0}Y(A,z)z^ndz
\end{equation}

In the affine Kac-Moody case, given a $\widehat{\mathfrak{g}}$-module $M_x$ we define its space of co-invariants as the quotient $M_x/\mathfrak{g}_{out}(x).M_x$. The space of conformal blocks is the dual space 

\begin{equation}
(M_x/\mathfrak{g}_{out}(x).M_x)^{\vee}=Hom_{\mathfrak{g}_{out}(x)}(M_x,\mathbb{C})
\end{equation}

\noindent
The basis elements of $\widehat{\mathfrak{g}}$ are given by $J^a(z)$, and $J^a(z)dz$ is naturally a one form. Then $\phi \in M_x^*$ is a conformal block if and only if $\langle \phi, J^a(z) .A \rangle$ has regular singularity at $x$. 

If we are given a vertex algebra $V$ and a Harish-Chandra pair $(\mathcal{B}, \mathcal{G}_+)$ of what is called internal symmetries (in case of Hodge structure come from Mumford-Tate group), then to any such data one can attach a twisted space of co-invariants $H_{\tau}(X)$.  Here the Lie algebra action of $\mathcal{B}$ is induced by the Fourier coefficients of vertex operators in $V$. And $\mathcal{G}_+$ consists of certain symmetries of a geometric data $\mathcal{H}$ (which in our case is the local system of Hodge structure) on the punctured disc $D^{\times}$. As $\tau$ varies $H_{\tau}(X)$ combine into a twisted $D$-module $\Delta(V)$ on the moduli $\mathcal{M}_{\Phi}$ of the data $\Phi$ on $X$.

\end{itemize}

\end{document}